\documentclass{article}
\usepackage{arxiv}
\usepackage[utf8]{inputenc}
\usepackage[T1]{fontenc}
\usepackage[ukrainian,russian,english]{babel}
\usepackage[unicode, pdftex]{hyperref}
\usepackage{url}            % simple URL typesetting
\usepackage{booktabs}       % professional-quality tables
\usepackage{amsfonts}       % blackboard math symbols
\usepackage{nicefrac}       % compact symbols for 1/2, etc.
\usepackage{microtype}      % microtypography
\usepackage{lipsum}		% Can be removed after putting your text content
\usepackage{graphicx}
\usepackage{amsmath}
\usepackage{amsfonts}
\usepackage{amssymb}  
\usepackage{natbib}
\usepackage{doi}
\usepackage{color}

\title{Algorithms and topological invariants for dynamic systems. III. Algorithms for Recognition and Classification of 2-Dimensional Surfaces}

\author{ \href{https://orcid.org/0000-0002-7164-807X}{\includegraphics[scale=0.06]{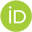}
\hspace{1mm}Alexandr O.~Prishlyak}\thanks{https://sites.google.com/view/prof-prishlyak, https://orcid.org/0000-0002-7164-807X} \\
	Department of Computer Methods of \\ 
	Mechanics and Control Processes\\
	Taras Shevchenko National University of Kyiv\\
	Kyiv, Ukraine \\
	\texttt{prishlyak@knu.ua} }

\renewcommand{\shorttitle}{Algorithms and topological invariants for dynamic systems. III. Surface Algorithms }

\newtheorem{theorem}{Theorem}

\newtheorem{problem}{Problem}

\hypersetup{
pdftitle={Algorithms and topological invariants for dynamic systems. II. Discrete Structures},
pdfauthor={Alexandr O. Prishlyak},
pdfkeywords={Topological classification, Morse function, Morse-Smale flow},
}

\begin{document}
\maketitle

\begin{abstract}
We construct algorithms and topological invariants that allow us to distinguish the topological type of a surface, as well as functions and vector fields for their topological equivalence.
In the first part (arXiv:2501.15657), we discused  basic concepts of diferential topology.   
In the second part (arXiv:2502.00506  ) we discused the main discrete topological structures used in the topological theory of dynamic systems.
In third part we construct algorithms that allow us  recognise  2-manifolds and 3-manifolds in simplicial complexes and regular CW-complex and detreminate topological type for 2-manifolds (connectivity, type of orientability, genus, number of boundary component).

\end{abstract}

% keywords can be removed

\section*{Introduction}

Topological properties of manifolds, functions, and dynamical systems are often studied through the construction of topological invariants that have a discrete nature, meaning they can be described using a finite set of integers, which allows for the use of computational techniques when working with them. Morse theory enables the construction of a cell complex structure on a manifold; however, continuous mappings used to attach cells are generally difficult to encode. Therefore, triangulations are more commonly used in computer modeling. Their drawback is the large number of simplices required to construct such structures. A compromise solution between these two structures is regular cell complexes. To find all possible structures under study, efficient algorithms for their recognition are necessary. Topological invariants will be useful if there are efficient algorithms for their computation and comparison. The issues mentioned are addressed in low dimensions (up to four) in this tutorial.

The first part  (\cite{prish2025atids1}) discusses the main structures used in the topology of manifolds: vector fields, dynamical systems, Morse functions, cell decompositions, and the fundamental group.

The second part (\cite{prish2025atids2}) examines discrete structures for which computers can be used for their specification and manipulation.

In this paper, we construct algorithms to determine the topological properties of simplicial complexes and regular cell complexes.

In the first section, we describe an algorithm that can check a 2-dimensional complex for connectivity and find the number of its components. In the second section, we develop an algorithm for checking local planarity. This consists of checking local planarity on edges and checking local planarity at vertices. We also determine the number of boundary components of the surface. The third section describes an algorithm for checking the orientability of the surface. The fourth section presents the Euler characteristic, the genus of the surface, and defines its topological type. The fifth section discusses possible generalizations of the constructed algorithms to higher dimensions. The sixth section describes topological invariants of stratified 2-dimensional sets, where 2-dimensional cells are replaced by 2-strata, i.e., surfaces with boundaries. The seventh section introduces a rotation system for a graph embedded in a surface. A special case of such systems (chord diagrams) is discussed in the eighth section. The ninth section provides a list of regular cellular complexes for 2-dimensional surfaces. The tenth section presents a list of rotation systems and chord diagrams with a small number of edges. The eleventh section lists graphs embedded in surfaces with no more than three edges, as well as graphs embedded in spheres with four edges.

\section{Algorithm for Checking Connectivity}
For cell (simplicial) complexes, connectivity is equivalent to path-connectedness. Since a disk and a simplex are convex sets, any two points within them can be connected by a path (line segment) that belongs to that disk (simplex). Therefore, for simplicial and regular cell spaces, the existence of a path between any two points is equivalent to the existence of a path between simplices (cells).

\subsection{Algorithm for Checking the Connectivity of a Simplicial Complex and Finding Connected Components.}

Let the simplicial complex \( K \) is given as a list of simplices (full or reduced).

1. We  find the lists \( S_i \) of all vertices belonging to each of the connected components. Consider the first simplex from \( K \). We add all its vertices to the list \( S_1 \) and cross this simplex out from the list of simplices in \( K \).

2. Let \( v_0 \) be the first vertex in the list \( S_1 \). We  find all simplices in \( K \) that contain \( v_0 \). We add their vertices to \( S_1 \) and remove these simplices from the list \( K \).

3. For the next vertex in the list \( S_1 \), we  perform the procedure described in step 2 (adding vertices to the list \( S_1 \) and crossing out simplices from \( K \)). We then move on to the next vertex...

This process  concludes in one of two possible ways:

a) All simplices in the list \( K \) are crossed out, which is equivalent to \( S_1 \) containing all vertices of the simplicial complex. In this case, the simplicial complex \( K \) is connected;

b) For the last vertex in the list \( S_1 \), there are no uncrossed simplices in \( K \), while there are still uncrossed simplices in \( K \). Hence, the simplicial complex is disconnected. To find the connected components, we take the following step.

4. In the list \( S_2 \), we add the vertices of the first uncrossed simplex from \( K \) and repeat steps 2 and 3. If the process ends again as in b), we  add the subsequent vertices to list \( S_3 \) and continue to add vertices to this and subsequent lists until there are no more vertices left. As a result, the number of non-empty \( S_i \) equals the number of connected components.

\begin{theorem} \cite{hat02, prish-modtop06} A cellular complex is connected if and only if its 1-skeleton is connected.
\end{theorem}

\textbf{Consequence.} The algorithm for checking connectivity and finding components can be simplified by considering its 1-skeleton instead of the complex. Since the 1-skeleton is a (topological) graph, the described algorithm is equivalent to the algorithm for finding connected components of a graph.

We demonstrate the operation of the algorithm with an example.

Let $$K=\{ \{1\},\{2\},\{3\},\{4\},\{5\},\{6\},\{1,2\},\{3,5\},\{2,4\},\{1,4\} \}.$$

1. $S_1=\{1\}$ (we add the vertices of the first simplex -- 1 to $S_1$);

2. $S_1=\{1,2,4\}$ (we add the vertices of the simplices containing 1 to $S_1$. Since there are no other simplices containing 1, 2, and 4, we move on to $S_2$);

3. $S_1=\{1,2,4\}, S_2=\{3\}$ (we add the vertices of the first simplex that does not contain 1, 2, or 4 to $S_2$);

4. $S_1=\{1,2,4\}, S_2=\{3, 5\}$ (we add the vertices of the simplices that do not contain 3 to $S_2$ and move on to $S_3$);

5. $S_1=\{1,2,4\}, S_2=\{3, 5\}, S_3=\{6\}$. Result: 3 connected components.

For a regular cell complex, the 1-skeleton is a graph (and thus a simplicial 1-complex); therefore, according to the last theorem, we can apply the algorithm for simplicial complexes to determine its connectivity.

\begin{problem}
Demonstrate the operation of the connectivity checking algorithm for the following spaces:

1) torus; 2) torus with a hole; 3) torus with two holes; 4) sphere with three holes; 5) Klein bottle; 6) M\"obius strip with a hole; 7) Klein bottle with two holes.
\end{problem}

\section{Algorithm for checking 2-dimensional complex for local planarity}
\textbf{Definition.} %/\begin{definition} \index{locally planar}
A topological space $X$ is called \textit{locally planar} if for every point $p \in X$, there exists a neighborhood $U$ that is homeomorphic to an open set of the plane or the upper half-plane.
%\end{definition}

The algorithm for checking the local planarity of a simplicial complex consists of two parts: 1) edge check; 2) vertex check.

1) \textbf{Edge check}: For each edge in the edge list, we find all 2-simplices containing it. There are three possible scenarios: a) the edge is contained in two 2-simplices, in which case it is internal; b) the edge is contained in one 2-simplex, in which case it is boundary; c) the edge is not contained in any 2-simplices or is contained in more than two 2-simplices, in which case the complex is not locally planar. If all edges are internal or boundary, then the complex is locally planar at the edges.

2)\textbf{ Vertex check} (for complexes that are locally planar at the edges). We fix a vertex. Let this vertex be $v_0$. We consider all 2-simplices that contain vertex $v_0$. For each of them, we find the edge opposite to $v_0$ (removing it from the trio of vertices including $v_0$). If the resulting list of edges, together with their endpoints, forms a connected complex (we apply a connectivity check algorithm), then the complex is locally planar at this vertex; otherwise, the complex is not locally planar. The complex is locally planar if it is locally planar at all edges and vertices.

Note that the connectivity of the set of opposite edges is equivalent to the ability to form a path containing all vertices involved. To do this, we take the first edge, let it be $\{v_1,v_2\}$. We remove it from the list of opposite edges. We then look for other edges in the list that contain $v_2$. If such an edge is unique, for example $\{v_2,v_3\}$, we append vertex $v_3$ to the sequence $v_1 \to v_2$ to the right of $v_2$ (resulting in $v_1 \to v_2 \to v_3$). We remove the edge $\{v_2,v_3\}$ from the list. If there are several edges containing $v_2$, then local Euclidean property is violated at edge $\{v_0,v_2\}$, which is impossible. If such an edge does not exist, we move to vertex $v_1$. We find another edge that contains it and its second vertex (let's call it $v_4$) and write it on the left.

We obtain, $v_4\to v_1\to v_2$, and we cross out this edge from the list. In the next step, we perform similar procedures with the right end of the sequence of vertices; if there is no second edge for it, we do so with the left end. After adding each subsequent vertex to the sequence, we check whether the left end of the sequence equals the right. The process of writing out the sequence ends either if there are no edges continuing it to the right or left, or if the left end of the sequence equals the right end. If the resulting sequence contains all vertices adjacent to $v_0$ (this is equivalent to saying that all edges are crossed out from the list of opposite edges), then the complex is locally planar at the vertex $v_0$. Otherwise, it is not locally planar.

\textit{For a regular cell complex}, we have the same two checks: 1) the check on edges is the same (simplices are replaced by cells), and 2) the check at vertices is analogous (instead of the opposite edge, we take the list of all 1-cells in the boundary of the 2-cell, whose endpoint is not the vertex $v_0$).

Note that boundary edges (and their endpoints) constitute the boundary of the surface. To find the number of components of the boundary, we apply a connectivity check algorithm.

\textbf{Examples.}

1. Let us consider the simplicial complex $$K=\{\{0\},\{1\},\{2\},\{3\},\{0,1\},\{0,2\},\{0,3\},\{1,2\},\{1,3\},$$ $$\{0,1,2\},\{0,1,3\}\}.$$

1) Local planarity on edges: $\{0,1\}$ is an interior edge (it is contained in each of the two 2-simplices), while the rest of the edges are boundary edges (they are contained only in one 2-simplex). Therefore, the complex is locally planar with respect to the edges.

2) Local planarity at vertices: a) for vertex $\{0\}$, the opposite edges $\{1,2\},\{1,3\}$ share a common vertex $\{1\}$, thus they form a connected set. The path that passes through all vertices is $3 \to 1 \to 2$; b) for vertex $\{1\}$, the corresponding path is $3 \to 0 \to 2$; c) for $\{2\}$ -- the path is $0 \to 1$; d) for $\{3\}$ -- also the path $0 \to 1$.

Therefore, the complex is locally planar with the boundary $[0, 2, 1, 3]$ (composed of boundary edges: $\{0,2\},\{2,1\},\{1,3\},\{3,0\}$).

\begin{problem}
Demonstrate the operation of the local planarity checking algorithm for the following spaces:

1) torus; 2) torus with a hole; 3) torus with two holes; 4) sphere with three holes; 5) Klein bottle; 6) M\"obius strip with a hole; 7) Klein bottle with two holes.
\end{problem}

\section{Algorithm for Checking the Orientability of Surfaces}

In this section, we  consider connected compact surfaces. They are defined using connected simplicial or regular cell complexes that are locally planar. Surfaces are divided into two types: orientable and non-orientable. Orientable surfaces are those for which all 2-cells (2-simplices) can be oriented such that any two adjacent 2-cells (neighboring 2-simplices) have consistent orientations. This means that the orientations they induce on the common edge (1-cell) are opposite. Recall that the orientation of a simplex can be given by an ordering of its vertices. The orientation of the 2-simplex $\{v_0,v_1,v_2\}$ induces the following orientations on its edges: $\{v_0,v_1\},\{v_1,v_2\},\{v_2,v_0\}$. For a two-dimensional cell, its orientation is determined by the orientation (direction of traversal) of its boundary. In this case, the restriction on each 1-cell from the boundary defines the orientation of that 1-cell.
From this definition, an algorithm for checking orientability follows.

We fix the orientation of the first 2-cell (2-simplex). For the second cell, which shares a common edge with the first cell (i.e., a one-dimensional cell that lies on the boundary of both the first and the second cell), we  choose an orientation that aligns with the first cell. For all subsequent cells that share at least one common edge with any of the examined cells, we choose an orientation consistent with the orientation of the first cell in the list that it intersects with and check whether this orientation is consistent with the orientations of the previous cells. If it is consistent, we move to the next cell; if not, the process stops with result ''The surface is non-orientable''. If we reach the last cell and its orientation is consistent with the previous ones (and thus all cell orientations are aligned), then the complex is orientable.

%\begin{example}
\textbf{Example.} The simplicial complex $$K=\left\langle \{ 0,1,2 \}, \{ 0,1,3 \}, \{ 0,2,3 \}, \{ 1,2,3 \} \right\rangle.$$

1) We align the orientation of the second simplex with the first on the common edge $\{0,1\}$, resulting in: $\{ 0,1,2 \}, \{ 1,0,3 \}$;

2) We align the orientation of the third simplex with the first on the edge $\{0,2\}$, resulting in: $\{ 0,1,2 \}, \{ 1,0,3 \}, \{ 0,2,3 \}$;

3) We check the consistency of the orientations of the second and third simplices on the common edge $\{0,3\}$: they are consistent;

4) We align the orientation of the fourth simplex with the first on the edge $\{1,2\}$, resulting in: $\{ 0,1,2 \}, \{ 1,0,3 \}, \{ 0,2,3 \}, \{ 2,1,3 \}$;

5) We verify the consistency of the orientation of the fourth simplex with the orientations of the second and third simplices on the common edges $\{1,3\}$ and $\{2,3\}$, respectively -- they are consistent.

Therefore, the simplicial complex $K$ is oriented. One of the two orientations is: $\{ 0,1,2 \}, \{ 1,0,3 \}, \{ 0,2,3 \}, \{ 2,1,3 \}$.
%\end{example}

%\begin{example}
\textbf{Example.}
The cell complex $$RCC=\left\langle [ 0,1,2, 3] , [ 2,3,4,5], [0,1,4,5] \right\rangle$$

1) We align the orientations of the first and second cells on the common edge $\{2,3\}$, obtaining $[ 0,1,2, 3] , [ 5,4,3,2]$;

2) We align the orientations of the first and third cells on the common edge $\{0,1\}$, obtaining $[ 0,1,2, 3] , [ 5,4,3,2] , [5,4,0,1]$;

3) We check the consistency of the orientations of the second and third cells on the common edge $\{4,5\}$ -- they are not consistent. Therefore, the complex $RCC$ is not oriented.
%\end{example}

\begin{problem}
Demonstrate  the orientability-checking algorithm for the following spaces:

1) torus; 2) torus with a hole; 3) torus with two holes; 4) sphere with three holes; 5) Klein bottle; 6) M\"obius strip with a hole; 7) Klein bottle with two holes.
\end{problem}

\section{Determination of the genus and topological type of a surface}
%\begin{definition} \index{genus of a surface}
\textbf{Definition.} The genus of a compact connected surface is defined as the maximum number of closed curves without intersection points, the complement of which is connected.
%\end{definition}

According to Jordan's theorem, any closed curve on a sphere divides it into two parts. Therefore, the sphere has genus 0. Similarly, a sphere with holes also has genus 0. For a torus, parallels and meridians (or any other curve that is not the boundary of a simply connected region on the torus) do not divide the torus. Cutting the torus along this curve yields a surface homeomorphic to a sphere with two holes, on which there are no closed curves that do not divide the surface. Thus, the genus of the torus is 1. If a M\"obius strip is cut along its central line, it results in a connected surface homeomorphic to a cylinder (a sphere with two holes). The genus of the M\"obius strip is 1. The M\"obius strip can be viewed as a projective plane with a hole, so the genus of the projective plane is also 1.
If a Klein bottle is represented as a sphere with two holes, each sealed with a M\"obius strip, we can cut them along the central lines, keeping the surface connected. Therefore, the genus of the Klein bottle is 2.

For any closed orientable surface $F_g$ of genus $g$, there is a formula that relates the genus to the Euler characteristic $\chi (F_g)$:

$$ g= \frac{2-\chi (F_g)}{2}.$$

For any closed non-orientable surface $N_g$ of genus $g$:
$$ g= 2-\chi (N_g).$$

If the surface has a boundary, then sealing each boundary component with a 2-disk (adding a 2-cell) increases the Euler characteristic by the number of boundary components, resulting in a closed surface. Thus, for an orientable surface $F_{g,b}$ of genus $g$ with $b$ boundary components, the following formula holds:

$$ g= \frac{2-\chi (F_{g,b})-b}{2}.$$

For a non-orientable surface $N_g$ of genus $g$ with $b$ boundary components:
$$ g= 2-\chi (N_{g,b})-b.$$

Therefore, the algorithm for checking whether a finite two-dimensional complex is a surface and determining its type is as follows:

1. We identify the components of connectivity and carry out the subsequent steps for each of them separately.

2. We check for local planarity at the vertices and edges and find the number of boundary components of the surface.

3. We determine the type of orientability.

4. We find the Euler characteristic and the genus.

The triplet consisting of the type of orientability, genus, and the number of boundary components is a complete topological invariant of a compact surface (defining the topological type of the surface).

\begin{theorem} A two-dimensional complex defines a surface homeomorphic to a sphere if and only if it is connected, locally planar, and has an Euler characteristic of 2.
\end{theorem}
%\begin{proof} 
\textbf{Proof.} The higher the genus, the lower the Euler characteristic. Similarly, each hole decreases the Euler characteristic by 1. Therefore, the maximum Euler characteristic for a surface of genus 2 without holes, that is, the sphere, is 2.
%\end{proof}

\textbf{Corollary.} A two-dimensional complex defines a surface homeomorphic to a 2-disk if and only if it is connected, locally planar, has a non-empty boundary, and an Euler characteristic of 1.

\textbf{Proof.} 
%\begin{proof} 
The 2-disk has the maximum Euler characteristic among surfaces with a non-empty boundary, which equals 1.
%\end{proof}
\begin{problem}
Determine the genus and topological type using the constructed algorithms for the following spaces:

1) torus; 2) torus with a hole; 3) torus with two holes; 4) sphere with three holes; 5) Klein bottle; 6) M?bius strip with a hole; 7) Klein bottle with two holes.
\end{problem}

\section{Algorithms for Higher-Dimensional Manifolds}

The algorithm for checking connectivity is the same for all dimensions. Therefore, let's move on to algorithms for checking local Euclidean properties (in dimension 2 -- planarity). It suffices to check local Euclidean properties at the vertices since if it is violated on the edges or simplices (cells) of higher dimensions, it will also be violated at the vertices lying on the boundary of that edge (simplex, cell). We consider the vertex $v_0$ of the simplicial complex and all the simplices containing this vertex and find the opposite faces to $v_0$ (removing $v_0$ from the vertex list of the simplex). The neighborhood of $v_0$ can be seen as a cone over the formed simplicial complex $S$. The complex is locally Euclidean at the vertex if $|S|$ is homeomorphic to a sphere or a disk of codimension 1. For a three-dimensional complex, the corresponding sphere and disk have dimension 2, and we can utilize the latest theorem and corollary for their verification. In higher dimensions (4 and above), the verification algorithm is a complex problem.

The algorithm for determining the boundary of a three-dimensional manifold is analogous to that for surfaces: the boundary consists of those 2-simplices that are faces of only one 3-simplex.

The algorithm for determining orientability is also similar to that for surfaces: we orient the first 3-simplex, then we align the orientations of the 3-simplices that share a common 2-face with this simplex, subsequently aligning the orientations of the simplices that share a boundary with already considered simplices and checking whether all constructed orientations are consistent with one another.

Other methods for defining 3-manifolds and 4-manifolds using Heegaard diagrams and Kirby diagrams is described in subsequent sections.

\begin{problem}
Construct a regular cellular decomposition and demonstrate the operation of topological property checking algorithms for the following manifolds:

1) $S^3$; 2) $T^3$; 3) $S^1 \times S^2$; 4) the full pretzel of genus 2; 5) the full torus; 6) $S^2 \times [0,1]$; 7) $T^2 \times [0,1]$; 8) $\mathbb{R} P^2 \times [0,1]$.
\end{problem}

\section{SLW-graphs of stratified sets and graph embeddings in surfaces}

The construction of a regular two-dimensional cellular complex can be generalized to the case of stratified two-dimensional sets.

From now on, we understand a graph as a one-dimensional cellular complex, that is, a multigraph with multiple edges and loops.

%\begin{definition} \index{semigraph} 
\textbf{Definition.} A semigraph is defined as a disconnected union of a graph with circles. A semigraph is finite if it has a finite number of vertices, edges, and circles.

%\end{definition}
Note that a semigraph can be made into a graph by adding a vertex to each circle and transforming that circle into a loop.

%\begin{definition} \index{combinatorial stratified set} 
\textbf{Definition.}
\textit{A combinatorial stratified two-dimensional set} is a topological space obtained from a finite semigraph through the attachment of a finite number of compact surfaces with boundary via continuous mappings of the boundaries, which are either local homeomorphisms, with the exception of possibly a finite number of points that map to the vertices of the graph, or constant mappings to the vertices. The surfaces that are attached are called two-dimensional strata, the edges of the graph are one-dimensional, and the vertices are zero-dimensional strata. The graph is called a 1-skeleton.
%\end{definition}

The concept of a combinatorial stratified set generalizes the concept of a combinatorial cellular complex, in which disks are replaced with manifolds with boundary. Examples of stratified two-dimensional sets include surfaces with finite embedded graphs.

\textbf{Definition.}
%\begin{definition} 
Two combinatorial stratified two-dimensional sets are called isomorphic if there exists a homeomorphism between them that preserves their stratification, that is, a homeomorphism whose restriction to 1-skeleta defines an isomorphism of graphs.
%\end{definition}

In this subsection, we address the problem of when an isomorphism of combinatorial stratified two-dimensional sets exists. The question of the existence of such an isomorphism is equivalent to the question of the possibility of extending an isomorphism of directed graphs to a homeomorphism of combinatorial stratified two-dimensional sets.

This problem arises in the topological classification of Morse-Smale vector fields on 2- and 3-dimensional manifolds.

It is also related to the possibility of extending isomorphisms of graphs to homeomorphisms of combinatorial stratified two-dimensional sets.

Let $G$ and $G'$ be directed graphs, and let $g: G \to G'$ be an isomorphism of graphs that maps the vertices $A_i$ of graph $G$ to the vertices $A'_i$, and the edges $B_j$ to $B'_j$. We denote $F_i$ (respectively, $F'_i$) as surfaces with boundaries that are glued to $G$ (and $G'$ respectively). Each circle that enters the boundary of surface $F_i$, corresponding to the vertices of the graph, is divided into arcs that map (homeomorphically) to the edges of graph $G$. In this case, the orientation of the edges of the graph determines the orientation of the arcs. We  assign the same letter to each arc as its corresponding edge in the graph.

On each surface $F_i$, we  fix an orientation if that surface is oriented. The orientation on each boundary component is assigned in accordance with the orientation of the surface in the case of oriented surfaces, and arbitrarily for non-oriented surfaces. For each circle on the boundary of the surface, traversing it according to the orientation, we will write a word consisting of letters $B_j^{\pm 1}$ corresponding to the edges encountered. A letter has a degree of $+1$ if the orientation of the corresponding arc matches the orientation of the circle, and $-1$ otherwise. Two words are called equivalent if one can be obtained from the other by a cyclic permutation of the letters. This corresponds to a different choice of the starting point for traversing the circle. Words are called reverse if one can be obtained from the other by writing the letters in reverse order and inverting the degrees, and possibly by cyclic permutation. This corresponds to traversing the circle against the orientation.

For each surface $F_i$, we  compile a list consisting of the following elements: 1) the number $n_i$ equal to the genus of the surface $F_i$ if the surface is oriented, and $-n_i$ otherwise; 2) the words written while traversing the boundary of the surface according to the orientation. We  call two such lists equivalent if the numbers $n_i$ match and there is a one-to-one correspondence between the letters of the words such that replacing the letters in the corresponding words yields equivalent or reverse words. Moreover, if the surface is oriented, then all words are either simultaneously equivalent or simultaneously reverse.

Thus, for a combinatorial stratified two-dimensional manifold, we have constructed a \textit{set of lists of words} (SLW) such that each list corresponds to one surface $F_i$. We  call two such sets equivalent if there exists a one-to-one mapping between their letters and a one-to-one correspondence between the letters and the lists such that the corresponding lists are equivalent under the given letter mapping.

%\begin{comment}
%\textbf{Examples.} Let us consider the graphs on the sphere as shown in Fig. 2.7 (here the sphere is represented as a plane with a point at infinity).

%The first SLW takes the form

%{{0, c}, {0, dbc-1b-1}, {0, aa-1d -1}},

%while for the second it is
%{{0, e}, {0, fh-1f -1e-1}, {0, hgg-1}}.

%Each list contains only one word, as the graphs are connected. We  apply the following letter substitutions to the first SLW:

%a ? g,
%b ? f -1,
5c ? e,
%d ? h-1.

%This results in the following SLW:

%{{0, e}, {0, h-1f -1e-1f}, {0, gg-1h}}.

%By cyclically permuting the letters in the second and third lists, we obtain the second SLW (SLW for G'). Thus, these two SLWs are equivalent.
%\end{comment}

\begin{theorem}(\cite{prishlyak1997graphs}) Let $G$ be a directed semi-graph that is a 1-skeleton of a combinatorial stratified two-dimensional set $K$, and let $G'$ be the 1-skeleton of $K'$, with $g: G \to G'$ being a graph isomorphism that maps the vertices $A_i$ of graph $G$ to the vertices $A'_i$, and the edges $B_j$ to $B'_j$. Then the graph isomorphism extends to a homeomorphism of stratified sets if and only if replacing $B_j$ with $B'_j$ in the list of words for $K$ results in a set equivalent to the list of words for $K'$.
\end{theorem}
\textbf{Proof.} \textit{\textit{Necessity.}} From the construction of the list of words, it follows that a homeomorphism of combinatorial stratified sets, when restricted to each surface $F_i$, defines an equivalence of the corresponding list of words.

\textit{Sufficiency.} The equivalence of lists allows us to construct homeomorphisms between the surfaces $F_i$ and $F_i'$, whose restrictions on the boundaries coincide with the isomorphism between the graphs. Indeed, these surfaces have the same genus and number of boundary components, and in the case of oriented surfaces, the orientation of the boundary is consistent with the orientation of the surface. We have that the graph isomorphism extends to each of the $F_i$, and thus to the entire set $K$.

A graph together with set of  list of words (SLW), where the letters correspond to the edges of the graph, is called an SLW-graph. Two SLW-graphs are called isomorphic if there exists a graph isomorphism such that replacing the letters of the second SLW with the letters corresponding to the isomorphism from the first SLW results in an SLW equivalent to the first SLW.

\textbf{Remark.} Each SLW-graph defines a combinatorial stratified two-dimensional set. Indeed, the list of words defines (up to homeomorphism) a compact surface with boundary and its attachment to the graph, which is a local homeomorphism except at the points corresponding to the vertices.

\subsection{SLW-graphs of Closed Surfaces}

Let us consider the question of when an SLW graph defines a surface. Let $K$ be a two-dimensional stratified set (complex) obtained by attaching surfaces to the graph corresponding to lists. Since each list of words defines a surface with a boundary, the structure of the surface (local homeomorphism to the plane) of the complex $K$ can only be violated at the gluing points, that is, at the edges and vertices of the graph.

1) The condition that the complex $K$ is locally planar at the interior points of the edges is equivalent to the fact that two edge parts are attached to each edge. This means that each letter or its inverse appears exactly twice in all words from the SLW.

2) Assuming the first condition holds, for each vertex, we consider the set of edges incident to it. We say that two edges are neighboring if they belong to the boundary of one of the attached surfaces. This condition is equivalent to the existence of a word where the letters corresponding to the edges, or their inverses (depending on the orientation of the edges and surfaces), are adjacent or are the first and last letters of the word. We say that two edges are equivalent if there is a chain of neighboring edges connecting them. Then, the condition for the complex $K$ to be locally planar at the vertices of the graph is equivalent to the requirement that for each vertex all edges incident to it are equivalent. If this condition is not met, the neighborhood of such a vertex is homeomorphic to a bouquet of planes, where the number of planes in the bouquet equals the number of equivalence classes.

\begin{problem}
Construct an SLW for graphs with three edges embedded in a sphere (see Fig. \ref{gs2-3}).
\end{problem}

\section{Rotation System of the Embedded Graph}

The rotation system is a dual construction for SLW-graphs of cellular embeddings of graphs (all surfaces $F_i$ are 2-disks).

For each vertex, we  fix the orientation in its regular neighborhood. If the surface is oriented, we  choose the orientation of the surface at the vertices. We  compile an ordered list of edges \( R \) by traversing the boundary of the vertex's neighborhood according to the orientation, starting from an arbitrary point (rotating around the vertex according to the orientation). For a non-oriented surface, we  also compile a consistency vector \( u \), consisting of + or - signs, corresponding to each edge and dependent on whether there is an orientation of the neighborhood edge that aligns with the orientations at its ends.

Analogous to the SLW, for a closed surface, each edge appears in all lists twice. If the list \( u \) contains all + signs, then the surface is oriented.

%\begin{example}
\textbf{Example.} 
\( R = \{ \{1,1,2\}, \{2,3,3\}\} \), \( u=\{ +, +, +\} \) defines a rotation system of a graph on a sphere, which has two vertices; edge 1 is a loop at the first vertex, edge 2 connects the vertices, and edge 3 is a loop at the second vertex.
%\end{example}

\begin{figure}[ht!]
\center{\includegraphics[width=0.6\linewidth ]{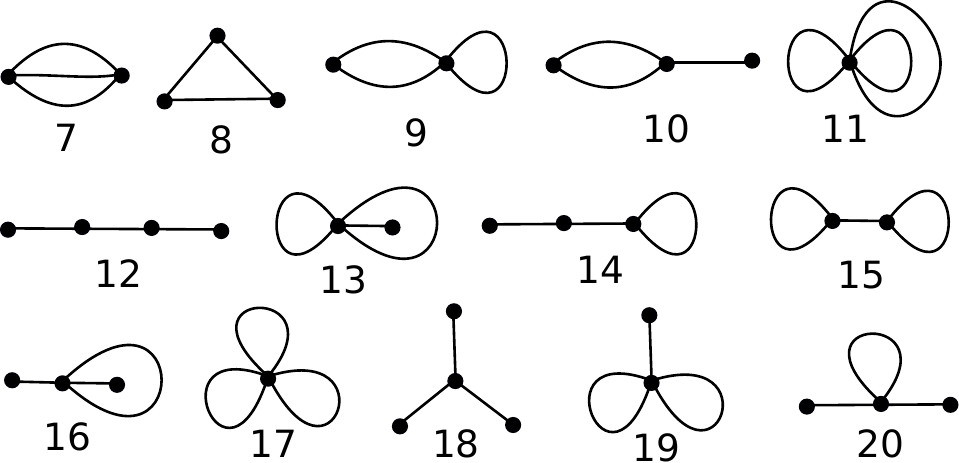}} 
\caption{Graphs with three edges on $S^2$
 }
\label{gs2-3}
\end{figure}

\textbf{Example.}
%\begin{example}
\( R = \{ \{1,1\}\} \), \( u=\{ - \} \) defines a rotation system of a graph on the projective plane, which has one vertex and one edge -- a loop, whose regular neighborhood is homeomorphic to a M\"obius strip.
%\end{example}
\begin{problem}
Create rotation systems for graphs with three edges, embedded in a sphere (see Fig. \ref{gs2-3}).
\end{problem}

\section{Chord Diagrams}
%\begin{definition}
\textbf{Definition.} A Hamiltonian cycle in a graph is defined as a cycle that passes through each vertex of the graph exactly once.
%\end{definition}

%\begin{definition} \index{chord diagram}
\textbf{Definition.} A chord diagram is a regular graph of degree 3  with a distinguished Hamiltonian cycle on it. Two diagrams are said to be isomorphic if there exists an isomorphism of graphs that maps the Hamiltonian cycle to a Hamiltonian cycle.
%\end{definition}

To represent chord diagrams, a Hamiltonian cycle is shown as a circle, and the edges not belonging to the Hamiltonian cycle are represented as chords.

We  number the chords  with the numbers 1, 2, 3, ... Starting from an arbitrary vertex, we  move along the cycle and write down an ordered set of numbers of the chords corresponding to the vertices we encounter. Each number  occur twice. The resulting set, which we  call the chord diagram code, uniquely defines the chord diagram up to isomorphism.

For example, if the chord diagram contains two intersecting chords, it is represented by the code $Ch=\{1,2,1,2\}$.

If a cell complex embedded in a closed oriented surface has one vertex, let us consider the regular neighborhood of the vertex. Its boundary is homeomorphic to a circle, which we  treat as a distinguished Hamiltonian cycle. The vertices of the chord diagram are the intersection points of the edges of the graph with this cycle. Two vertices are connected by a chord if they belong to the same edge. The resulting chord diagram has a code that coincides with the rotation system of the original graph.

Chord diagrams are widely used in modern low-dimensional topology, especially in knot theory.

Let $D$ be an arbitrary chord diagram with $2n$ points. We number them from 1 to $2n$ in a clockwise direction, starting from a certain fixed point on the double graph. We construct a permutation $\alpha$, obtained from the transposition product $(a_j, b_j), j = 1, \ldots, n$, where $a_j, b_j$ are the numbers of the points connected by a chord in the diagram $D$:

$$ \alpha = (a_1, b_1)(a_2, b_2) \ldots (a_n, b_n). $$

The choice of the starting point can lead to different permutations. For example,% the 4-diagram depicted in Figure 1 corresponds to the permutation
$\alpha_1 = (1,6)(2,8)(3,7)(4,5)$, and
$\alpha_2 = (1,7)(2,6)(3,4)(5,8)$ corresponds to the same chord diagram.

A \textit{nice} chord diagram is one in which all chords are based on a canonical circle (a single circle centered at the origin in $R^2$) with a fixed numbering of $2n$ points on it, which are the vertices of a regular polygon. For instance, we can number them counterclockwise, starting from the point (1,0). Then, each diagram $D(\alpha)$ uniquely determines the permutation $\alpha$ and vice versa.

\begin{problem}
Find all chord diagrams with 4 edges and for each of them find the inverse.
\end{problem}

\section{List of Regular Cell Complexes }
Let us provide examples regular cell complexes (RCC) for compact surfaces of low genus and a small number of boundary components.

Closed surfaces:

a) Sphere $S^2$:  $RCC=<012, 023, 031,132>$

b) Torus $T^2$: $RCC=<0153,0284,  0362,
 0471, 1265, 1782, 3486,$ $ 3574, 5687>$,

c) $\mathbb{R}P^2$:  $RCC=<013, 0154, 024, 0253, 125, 1243, 345>$

d) $Kl$: $RCC=<0153,0184,  0362,
 0472, 1265, 1278, 3486,$ $ 3574,$ $5687>$.

Compact surfaces with boundary:

a) $S^1\times [0,1]$: $RCC=<0154, 0231, 2453>,$

b) $Mo$: $RCC=<0132, 0154, 2354>,$

c) $F_{0,3}$: $RCC=<0123X5,0671, 2893, 5X9876>$,

d) $F_{1,1}$: $RCC=<01X523,0671, 2893, 5X9876>$.

\begin{figure}[ht!]
\centerline{\begin{minipage}[h]{0.2\linewidth}
\center{\includegraphics[width=1\linewidth ]{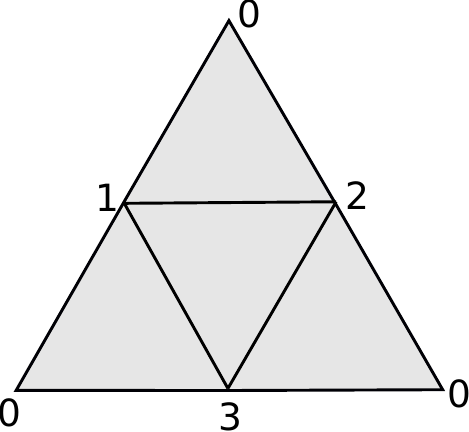}} a) \\
\end{minipage}
\ \ \
\begin{minipage}[h]{0.2\linewidth}
\center{\includegraphics[width=1\linewidth]{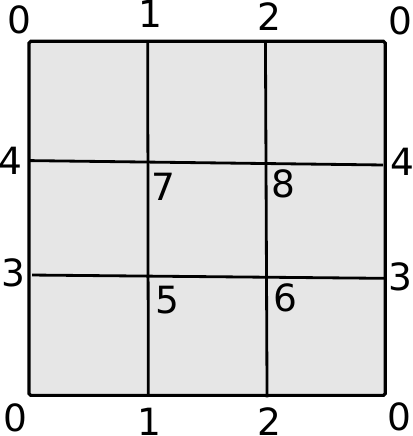}} \\b)
\end{minipage}
\ \ \ \
%}
%\vspace{0.4cm}
%\centerline{
\begin{minipage}[h]{0.24
\linewidth}
\center{\includegraphics[width=1\linewidth]{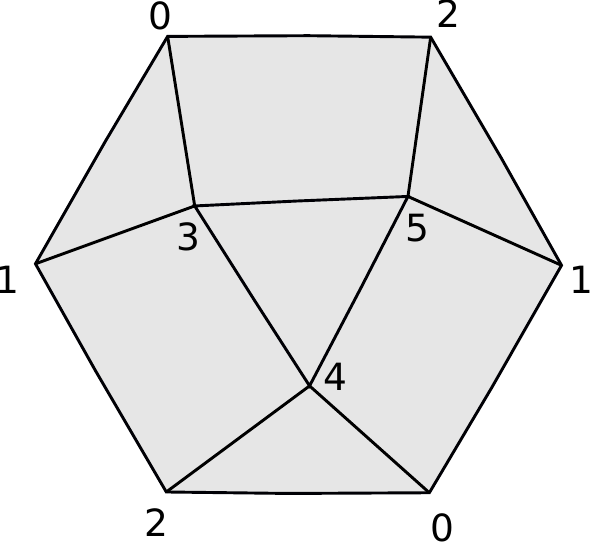}} c) \\
\end{minipage}
\ \ \ \
\begin{minipage}[h]{0.2\linewidth}
\center{\includegraphics[width=1\linewidth]{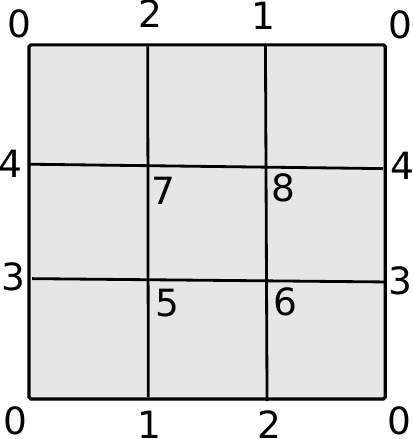}} d) \\
\end{minipage}}
\caption{Closed surfaces: a) sphere $S^2$, b) torus $T^2$, c) projective plane $\mathbb{R}P^2$, d) Klein bottle $Kl$.}
\label{pkk1}
\end{figure}

\newpage

\begin{figure}[ht!]
\begin{center}
\begin{minipage}[h]{0.25\linewidth}
\center{\includegraphics[width=0.9\linewidth ]{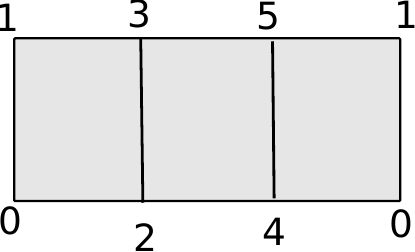}} \\ a)
\end{minipage}
\ \ \ \ \ \ \ \ 
\begin{minipage}[h]{0.25\linewidth}
\center{\includegraphics[width=0.9\linewidth]{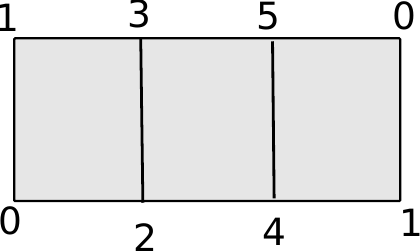}} \\b)
\end{minipage}
\ \ \ \ \ \ \ \ 
%\vfill
\begin{minipage}[h]{0.30\linewidth}
\center{\includegraphics[width=1\linewidth]{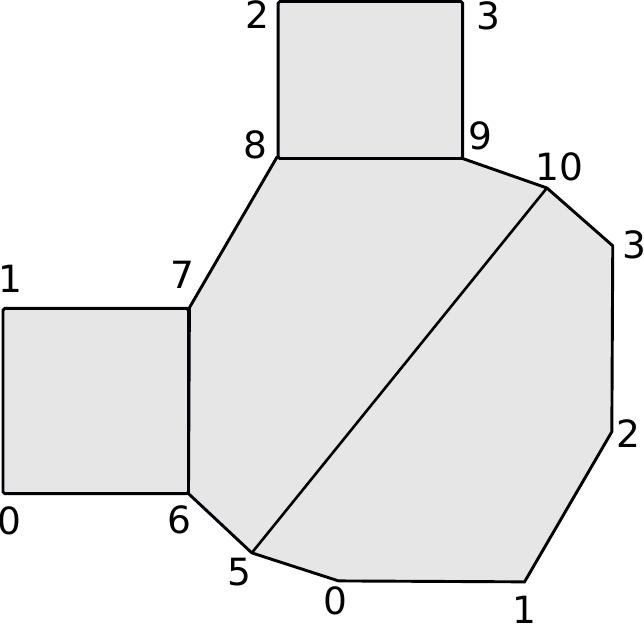}} c) \ \ \ \ \ \ \ \ 
\end{minipage}
\ \ \ \ \ \ \ \ \ 
\begin{minipage}[h]{0.30\linewidth}
\center{\includegraphics[width=1\linewidth]{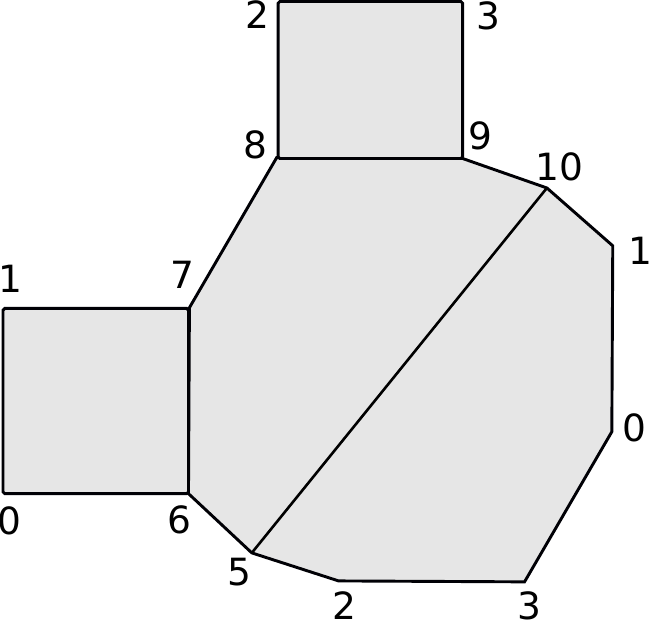}} d) \\
\end{minipage}
\ \ \ \ \ \ \ \ 
\begin{minipage}[h]{0.30\linewidth}
\center{\includegraphics[width=1\linewidth]{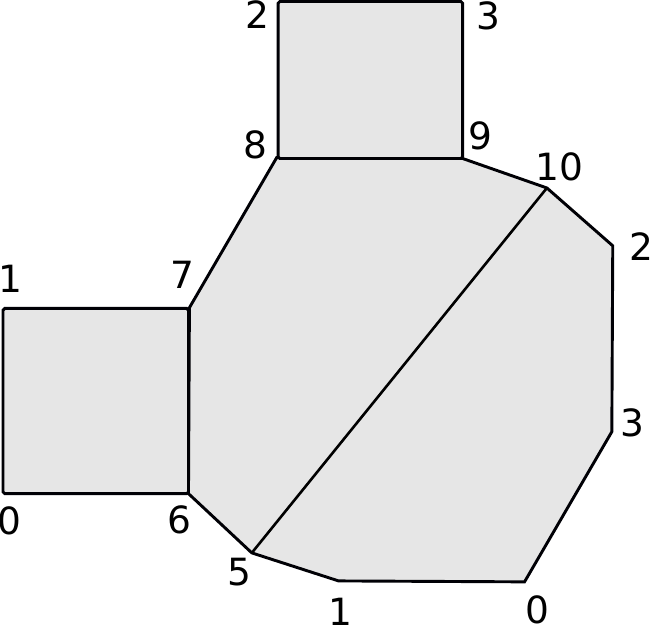}} e) \\
\end{minipage}
\ \ \ \ \ \ \ \ 
\begin{minipage}[h]{0.30\linewidth}
\center{\includegraphics[width=1\linewidth]{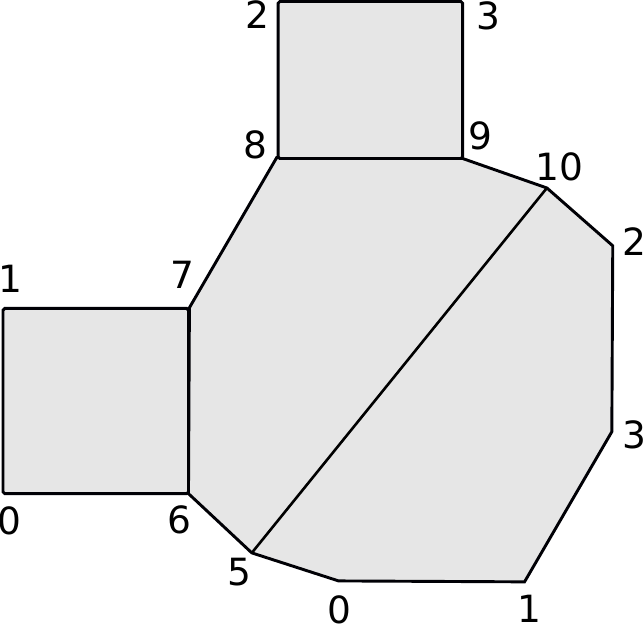}} f) \\
\end{minipage}
\end{center}
\caption{Surfaces with boundary: a) cylinder $S^1 \times [0,1]$, b) M\"obius strip $Mo$, c) sphere with three holes $F_{0,3}$, d) torus with a hole $F_{1,1}$, e) Klein bottle with a hole $N_{2,1}$, f) M\"obius strip with a hole $N_{1,2}$.
}
\label{pkk2}
\end{figure}

\begin{figure}[ht!]
\begin{minipage}[h]{0.45\linewidth}
\center{\includegraphics[width=0.80\linewidth ]{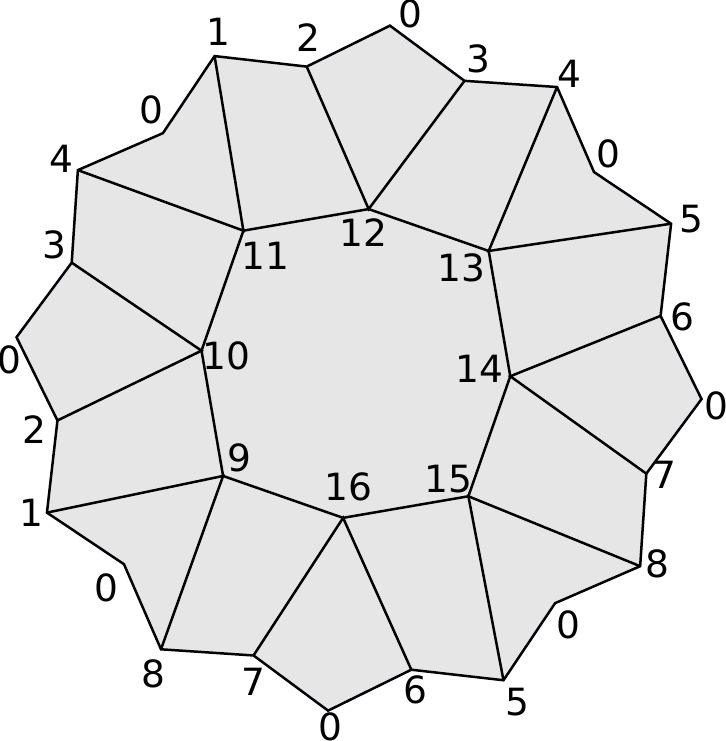}} \\ a) \\
\end{minipage}
\ \ \ \ \ \ \ \ 
\begin{minipage}[h]{0.45\linewidth}
\center{\includegraphics[width=0.7\linewidth]{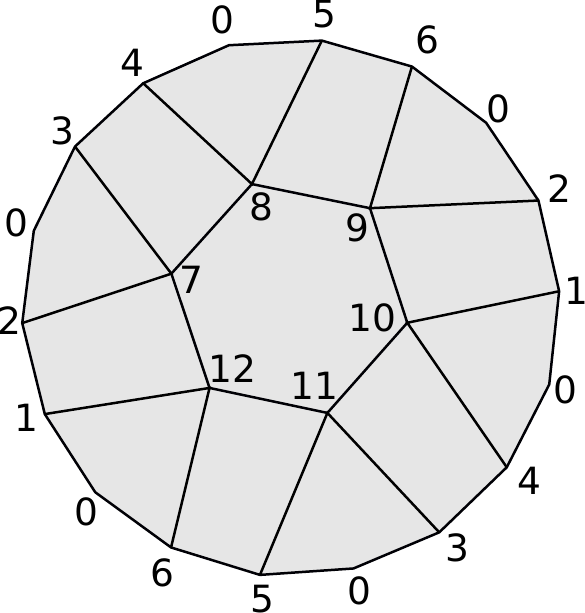}} \\ b) \\
\end{minipage}
\caption{Closed surfaces: a) orientable genus 2 surface $F_2$, b) non-orientable genus 2 surface 3 $N_3$.}
\label{pkk3}
\end{figure}

The RCC for an unoriented surface of genus 4 can be obtained from the RCC for $F_2$ by replacing the face with the boundary $[1,2,10,9]$ with a face having the boundary $[2,1,10,9]$, while leaving the other faces unchanged.

\newpage

\newpage

\begin{figure}[ht!]
\center{\includegraphics[width=0.4\linewidth]{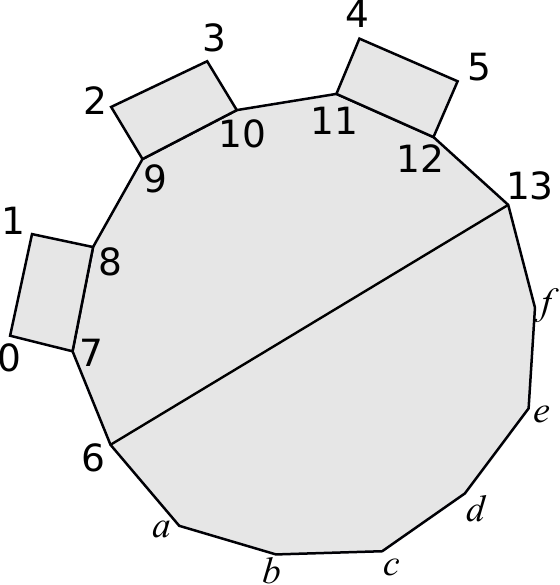}}
\caption{Surfaces with boundaries.}
\label{pkk4}
\end{figure}
Figure \ref{pkk4} shows such surfaces with boundaries:

1) a sphere with 4 holes $F_{0,4}$, $a=0, b=1, c=2, d=3, e=4, f=5;$

2) a torus with 2 holes $F_{1,2}$, $a=0, b=1, c=4, d=5, e=2, f=3;$

3) $Mo$ with 2 holes $N_{1,3}$, $a=0, b=1, c=2, d=3, e=5, f=4;$

4) $Kl$ with 2 holes $N_{2,2}$, $a=0, b=1, c=3, d=2, e=5, f=4;$

5) a genus 3 surface with a hole $N_{3,1}$, $a=1, b=0, c=3, d=2, e=5, f=4.$

For closed surfaces of higher genus, RCC can be constructed similarly to Figure \ref{pkk3}. For instance, for an unoriented surface of genus 5, such a decomposition  consists of one decagon and 20 quadrilaterals.

For compact surfaces with boundaries of higher genus, RCC can also be constructed similarly to Figure \ref{pkk4}. For example, the decomposition of a sphere with 5 holes (as well as a torus with 3 holes and an oriented surface of genus 3 with a hole) can be constructed with the following cells: two decagons and four quadrilaterals.

%\newpage

\section{List of Rotation Systems and Chord Diagrams}

Rotation systems with no more than 3 edges on the sphere $S^2$:

\begin{minipage}{0.4\textwidth}
  $R_1=Ch_1=\{11 \}$,

$R_2=\{1,1 \}$,

$R_3=\{12, 12 \}$,

$R_4=Ch_2=\{1122 \}$,

$R_5=\{1, 12, 2 \}$,

$R_6=\{1, 122 \}$,

$R_7=\{123, 132 \}$,

$R_8=\{12, 13, 23 \}$,

$R_9=\{1123, 23 \}$,

$R_{10}=\{12, 132, 3 \}$,
\end{minipage}
\hfill
\begin{minipage}{0.4\textwidth}
 $R_{11}=Ch_3=\{ 123321\}$,

$R_{12}=\{1, 12, 23, 3 \}$,

$R_{13}=\{ 11232, 3\}$,

$R_{14}=\{112, 23, 3 \}$,

$R_{15}=\{112, 233 \}$,

$R_{16}=\{ 1213, 2, 3 \}$,

$R_{17}=Ch_4=\{112233 \}$,

$R_{18}=\{123,1,2,3 \}$,

$R_{19}=\{11223,3\}$,

$R_{20}=\{1123, 2, 3 \}$.

\end{minipage}

Rotation systems with no more than 3 edges on the torus $T^2$:

\begin{minipage}{0.4\textwidth}
 $R_{21}=Ch_5=\{1212 \}$,

$R_{22}=Ch_6=\{123123 \}$,

$R_{23}=\{123, 123 \}$,

$R_{24}=Ch_7=\{123132 \}$, 
\end{minipage}
\hfill
\begin{minipage}{0.4\textwidth}
  $R_{25}=\{1213, 23 \}$,

$R_{26}=Ch_8=\{112323 \}$,

$R_{27}=\{12123, 3 \}$.
\end{minipage}

Chord diagrams (rotation systems) with 4 chords on the double torus  $T^2 \# T^2$:

\begin{minipage}{0.45\textwidth}
$R_{28}=Ch_9=\{12341234\}$,

$R_{29}=Ch_{10}=\{12312434\}$,
\end{minipage}
\hfill
\begin{minipage}{0.45\textwidth}
 $R_{30}=Ch_{11}=\{12132434\}$,

$R_{31}=Ch_{12}=\{12123434\}$.
\end{minipage}

 %\newpage
Chord diagrams with six edges on a surface of genus 3. 
\begin{figure}[ht!]
\center{\includegraphics[width=0.78\linewidth ]{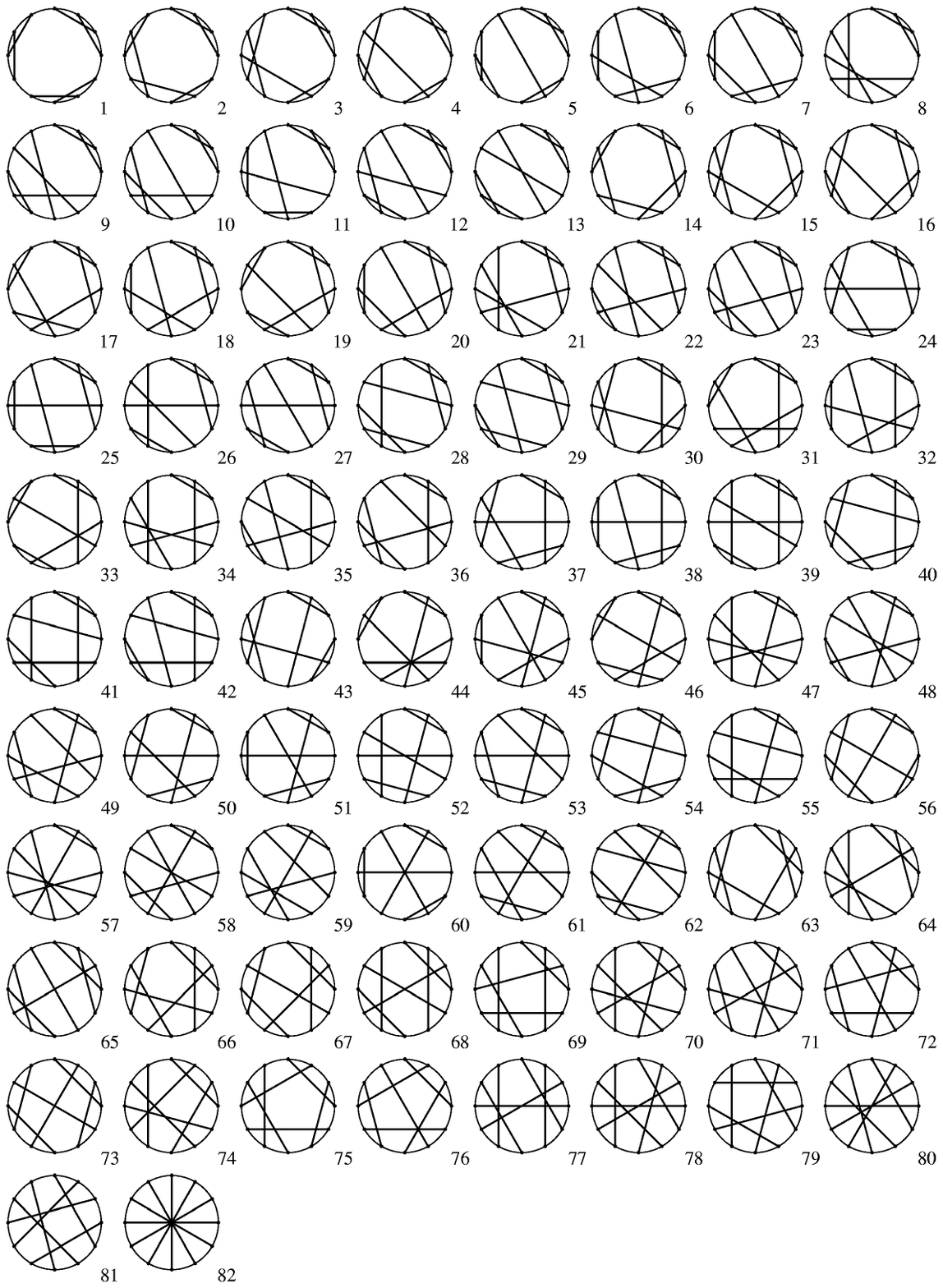}} 
\caption{Chord diagrams with six edges on   $F_3$
 }
\label{ch82}
\end{figure}

Rotation systems with no more than 2 edges on the projective plane
 $\mathbb{R} P^2$:

$R_{32}=Ch_{13}=\{\textcolor{blue}{11}\}, u=\{- \}$,

$R_{33}=Ch_{14}=\{\textcolor{blue}{1212}\}, u=\{- - \}$,

$R_{34}=\{\textcolor{blue}{1}2, \textcolor{blue}{1}2\}, u=\{- +\}$,

$R_{35}=Ch_{15}=\{11\textcolor{blue}{22}\}, u=\{+ - \}$,

$R_{36}=\{\textcolor{blue}{11}2,2\}, u=\{- +\}$.

Rotation systems (chord diagrams) with 2 edges on the Klein bottle
 $Kl$:

$R_{37}=Ch_{16}=\{\textcolor{blue}{1}2\textcolor{blue}{1}2\}, u=\{- + \}$,

$R_{38}=Ch_{17}=\{\textcolor{blue}{1122}\}, u=\{- - \} $.

\section{List of Embedded Graphs}
\begin{figure}[ht!]
\center{\includegraphics[width=0.6\linewidth ]{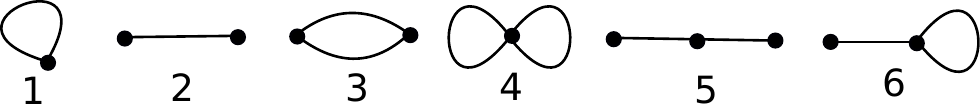}} 
\caption{Graphs with one and two edges on $S^2$
 }
\label{gs2-12}
\end{figure}

\begin{figure}[ht!]
\center{\includegraphics[width=0.65\linewidth ]{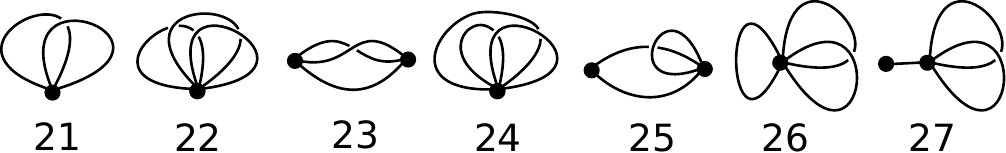}} 
\caption{Graphs with two and three edges on $T^2$
 }
\label{gt2}
\end{figure}

\begin{figure}[ht!]
\center{\includegraphics[width=0.65\linewidth ]{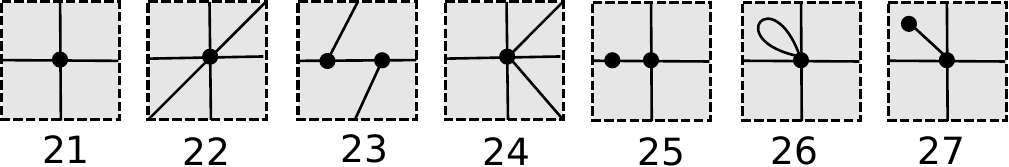}} 
\caption{Other representations of graphs with two and three edges on $T^2$
 }
\label{gt2b}
\end{figure}

\newpage

\begin{figure}[ht!]
\center{\includegraphics[width=0.65\linewidth ]{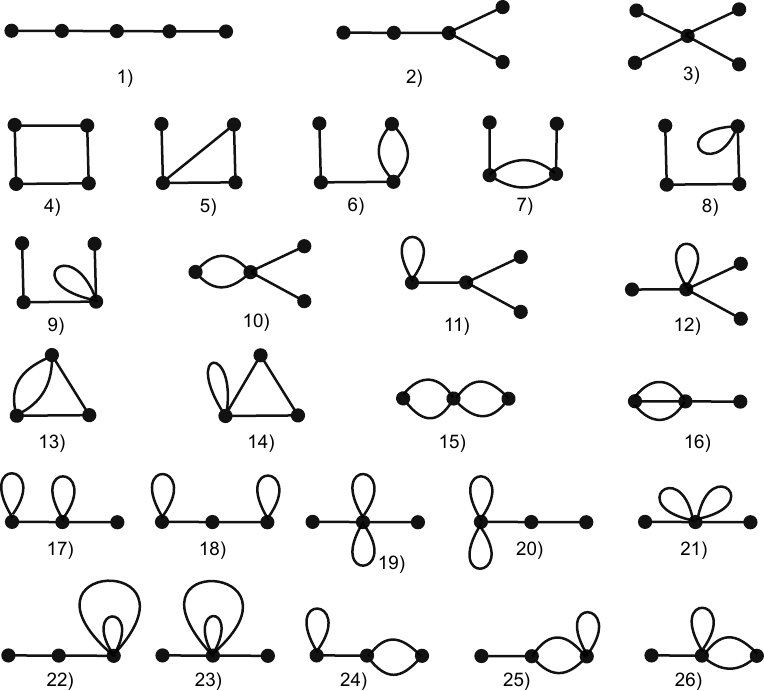}} 
\caption{Graphs with four edges on $S^2$
 }
\label{gs2-4}
\end{figure}

\section{Literature Review}

%\subsection*{Fundamentals of Differential Topology}

%A more detailed introduction to the main concepts of general (set-theoretical) topology can be found in \cite{bab15}.

%The concepts of smooth manifolds and smooth mappings, as well as vector fields, are covered in contemporary courses on differential geometry and topology \cite{bri94, lee00, bor95, hir79}.

%From Morse theory, one can recommend \cite{mil65, prish-tmorse02}.

%The concepts of cell complexes and fundamental groups are considered fundamental in algebraic (homotopy) topology \cite{hat02,  prish-modtop06}.%kos98,, spe71

%\subsection*{Discrete Topological Structures}

Simplicial complexes, Euler characteristics, and homology groups are fundamental concepts in homology theory courses \cite{hat02, prish-modtop06}. Regular cell complexes and discrete Morse functions are discussed in discrete Morse theory \cite{for98}. %, sco19}.

%\subsection*{Algorithms for Recognition and Classification}

%Algorithmic problems in the theory of low-dimensional manifolds were addressed in \cite{mf91}.

SLW-graphs were introduced in \cite{prishlyak1997graphs}.

Rotation systems are one of the core concepts in topological graph theory \cite{Peixoto1973, Gross1987,mohar2001graphs}.

Enumerating chord diagrams is covered in the works \cite{Kadubovski2021}.

%\subsection*{Topological Equivalence of Functions and Dynamical Systems}

%Classical textbooks on dynamical systems include \cite{Katok1995, Meiss2017, Palis1982, zol2023, kps2020}.

Topological properties of functions and vector fields have been studied in many works. Among them, we would like to highlight the works of the Kyiv topological school: 
textbooks \cite{prish2002Morse, prish-modtop06, prish2015top, prish23algkom}, % \cite{sha90}, 
scientific articles by the author on function topology \cite{Prishlyak1993, prishlyak1998, prishlyak1999equivalence, prishlyak2000conjugacy, prishlyak2001conjugacy, prish2001top, prishlyak2002topological1, prishlyak2002morse, prishlyak2003regular, prishter2024}, vector fields \cite{Prish1997vec, Prishlyak2001, prishlyak2002morse1, Prishlyak2002, prishlyak2003topological, prishlyak2003sum, prishlyak2005complete, Prishlyak2007}, and other geometric objects \cite{Prishlyak1994, prishlyak1997graphs, Prishlyak1999}, as well as the works of his students: K. Myshchenko \cite{prishlyak2007classification}; N. Budnytska \cite{Bud2008knu}; D. Lychak \cite{lychak2009morse}; A. Bondarenko \cite{Bond2012mfat}; O. Vyatychaninova \cite{VyatP2013Mol}; Bohdana Hladysh \cite{Hladysh2016, hladysh2017topology,% hladysh2019simple,
Hladysh2019, PLH2023}; A. Prus \cite{Prishlyak2017, prishlyak2020three, Prishlyak2021}; V. Kiosak \cite{KPL2022}; S. Bilun \cite{bilun2002closed, PBP23}%, Bilun2022, Bilun2022a}, 
V. Lisikevich \cite{LisP2013KNU}, I. Ivaniuk \cite{IvanPrish2014-func-deform3, IP2015knu},\cite{Ovtsynov2024}, D. Skotchko \cite{%Skochko2015, 
PS2016-F-atoms}, M. Loseva \cite{Losieva2017, Prishlyak2019, Prishlyak2020}, Z. Kybalko \cite{Kybalko2018}, K. Khatamian \cite{Hatamian2020}.%, L. di Beo \cite{Prishlyak2022}.

\bibliographystyle{unsrtnat}
\bibliography{prish}
\end{document}